\documentclass{amsart} \usepackage{amsmath,amssymb,amsfonts}
\usepackage[pdftex]{graphicx}

\usepackage{sty-ab,sty-va}

\renewcommand\eqref[1]{(\ref{#1})}

\date{August 15, 2012} 

\author{Valery Alexeev} \address{Department of Mathematics, University
  of Georgia, Athens, GA 30605, USA} \email{valery@math.uga.edu}
\author{Alexander Borisov} \address{Department of Mathematics, University of Pittsburgh,
301 Thackeray Hall, Pittsburgh, PA 15260, USA}
\email{borisov@pitt.edu} \title{On the Log Discrepancies in Toric Mori Contractions}

\begin{document}
\begin{abstract} It was conjectured by M\textsuperscript{c}Kernan and
  Shokurov that for all Mori contractions from $X$ to $Y$ of given
  dimensions, for any positive $\eps$ there is a positive $\delta$
  such that if $X$ is $\eps$-log terminal, then $Y$ is $\delta$-log
  terminal. We prove this conjecture in the toric case and discuss the
  dependence of $\delta$ on $\epsilon$, which seems mysterious.
\end{abstract}

\maketitle

\tableofcontents

\section{Introduction}
\label{sec:intro}

The main subject of this paper is the following 2003 conjecture of
James M\textsuperscript{c}Kernan:

\begin{Conjecture}[M\textsuperscript{c}Kernan]\label{conj:eps-delta}
  For fixed positive integers $m,n$ and a real number $\epsilon>0$
  there exists a $\delta=\delta_{m,n}(\epsilon)>0$ such that the
  following holds: Let $X$ be a $\bQ$-factorial variety, and $f\colon
  X\to Y$ be a Mori fiber space with $\dim Y=n$, $\dim X=m+n$.  Assume
  that $X$ is $\epsilon$-log terminal. Then $Y$ is $\delta$-log
  terminal.
\end{Conjecture}

A related stronger conjecture was suggested by V.V. Shokurov.  Let
$f\colon X\to Y$ be a proper surjective morphism with connected fibers
of normal varieties, so that $X/Y$ is of relative Fano type (see
definitions below) and let $\Delta$ be a $\mathbb Q$-divisor on $X$
such that $K_X+\Delta = f^*L$ for some $\bQ$-divisor $L$ on $Y$.

 By Kawamata's subadjunction formula
\cite{Kawamata_Subadjunction1, Kawamata_Subadjunction2}, see also
\cite{Ambro}, one has $K_X+\Delta = f^*(K_Y + R+ B)$, where $R$ is
the discriminant part, and $B$ is the ``moduli'' part, a $\mathbb Q$-divisor
defined only up to $\mathbb Q$-linear equivalence.

\begin{Conjecture}[Shokurov]\label{conj:Shokurov}
  In the above settings, assume that $(X,\Delta)$ is $\epsilon$-log
  terminal. Then there exists $\delta=\delta_{m,n}(\epsilon)>0$ and an
  effective moduli part $B$, such that $(Y,R+B)$ is
  $\delta$-log terminal.
\end{Conjecture}

Conjecture~\ref{conj:Shokurov} clearly implies
Conjecture~\ref{conj:eps-delta}: for a Mori fiber space consider a
large integer $N\gg0$ and a generic element $D$ of a very ample linear
system $-NK_X +f^*M$ for some $M$ on $Y$, and let $\Delta=\frac1{N}D$.
Then $K_X+\Delta = f^*L$ and for the minimal log discrepancies one has
\begin{displaymath}
  \mld (X, \Delta)  = \mld (X) + \frac1{N} 
  \quad \text{and} \quad
  \mld(Y,R+B) \le \mld(Y).
\end{displaymath}
Taking the limit $N\to\infty$ gives the implication.

We refer the reader to \cite{KollarMori_Book} for basic definitions
and results of the Minimal Model Program, some of which we briefly
recall below. 
For any normal variety $X$ for which
some positive multiple of the canonical class $K_X$ is Cartier,
one defines \emph{discrepancies} $a_i\in\bQ$ by the formula
\begin{displaymath}
  K_{X'} = \pi^* K_X + \sum a_i E_i, 
  \qquad \Exc(\pi) = \cup E_i,
\end{displaymath}
in which $\pi\colon X'\to X$ is a resolution of singularities, and
$E_i$ are the irreducible exceptional divisors of $\pi$.  The
\emph{log discrepancies} are the numbers $a_i^{\rm log}=a_i+1$. The
\emph{minimal log discrepancy} $\mld(X)$ is the infinum of log
discrepancies, going over all resolutions of singularities. Then
either $\mld(X)\ge 0$ or $\mld(X)=-\infty$. In the first case, variety
$X$ is log canonical, and $\mld(X)$ can be computed on any one 
resolution $\pi\colon X'\to X$ such that
$\Exc(\pi)$ is a normal crossing divisor.

A variety is said to be $\epsilon$-log terminal (abbreviated below to
$\epsilon$-lt) if its log discrepancies are $>\epsilon$, i.e. if for
ordinary discrepancies one has $a_i> -1+\epsilon$. Similarly, a
variety is $\epsilon$-log canonical if the log discrepancies are
$\ge\epsilon$.  In particular, $0$-log terminal is the same as 
Kawamata log terminal (klt), and $0$-log canonical is the same as log canonical.

We recall that $f\colon X\to Y$ is a \emph{Mori fiber space} if $f$ is
projective, $-K_X$ is $f$-ample, and the relative Picard number is
$\rho(X/Y)=1$. The assumption that $X$ is $\bQ$-factorial implies that
so is $Y$ (cf. \cite[Lemma 5-1-5]{KMM}).

Finally, a variety $X$ is called a variety of Fano type (FT) if there
exists an effective $\bQ$-divisor $D$ such that the pair $(X,D)$ is
klt and $-(K_X +D)$ is nef and big.

\medskip

There are numerous motivations for the above conjectures. The case
$\epsilon=\delta=~0$ of Conjecture~\ref{conj:eps-delta}, i.e. ``$X$ is
klt implies $Y$ is klt'' follows 
easily by cutting $X$ with $m$ general hyperplanes and reducing to a
finite surjective morphism. Even if $X$ is \emph{not} $\bQ$-factorial,
the implication ``$X$ is klt implies $(Y,\Delta)$ is klt for an
appropriate divisor $\Delta$'' is true, as proved
by Fujino \cite{Fujino}.

The first nontrivial case with $\epsilon>0$ appears when $\dim X=3$
and $\dim Y=2$, i.e. when $f\colon X\to Y$ is a singular conic
bundle. Mori and Prokhorov \cite{MoriProkhorov} considered the case
when $X$ is terminal. In this case, they proved
Iskovskikh conjecture which says that $Y$ must have at worst Du Val
singularities. This proves that one can take
$\delta_{1,2}(1)=1-c$ for any $c>0$.
Yuri Prokhorov also showed us several
examples of conic bundles of the form $(\bP^1\times \bA^2)/G \to
\bA^2/G$ for a cyclic group $G$ which indicate that
Conjecture~\ref{conj:eps-delta} is plausible.

Conjecture~\ref{conj:eps-delta} may also be viewed as the local
analogue of Borisov-Alexeev-Borisov (BAB) boundedness conjecture
\cite{BorisovBorisov,Alexeev_Boundedness} which says that for fixed
$n$ and $\epsilon>0$ the family of $n$-dimensional $\epsilon$-lt Fano
varieties is bounded.

Indeed, if $X$ happen to be Fano varieties, then the family of possible
$\epsilon$-lt varieties~$X$ is bounded by the BAB conjecture.  
Then
the family of possible varieties $Y$ must be bounded, so some
$\delta(\epsilon)>0$ must exist.
Vice versa, when trying to prove BAB conjecture by
induction, Conjecture~\ref{conj:eps-delta} naturally appears as one of the
steps. In this sense, it can be considered to be ``the local BAB
conjecture''.

\medskip

The main result of the present paper is the following

\begin{Theorem}\label{thm:main}
  Conjecture~\ref{conj:eps-delta} holds in the toric case, i.e. when
  $f\colon X\to Y$ is a morphism of toric varieties corresponding to a
  map of fans $(N_X,\Sigma_X)\to (N_Y,\Sigma_Y)$.
\end{Theorem}

Note that in the toric case, if one denotes by $\Delta$ the sum of
torus invariant divisors with coefficients 1, then 
 one has $K_X+\Delta = 0$ and the pair $(X,\Delta)$ is log canonical
 with $\mld(X,\Delta)=0$. Thus, the more general
 Conjecture~\ref{conj:Shokurov} does not fit the toric case very
 well. 

A very interesting question is to find the asymptotic of the function
$\delta(\epsilon)$ as $\epsilon\to 0$. Concerning this, we prove the
following:

\begin{Theorem}\label{thm:fiber-fixed}
  In the conditions of Theorem~\ref{thm:main}, suppose additionally
  that the generic fiber of $f$ is a finite, unramified 
in codimension one, toric quotient of a fixed
  toric Fano variety $P$.  
Then there exists a constant $C$ such that
  $\delta \ge C\cdot \eps^{m+1}$.
\end{Theorem}

On the other hand, we prove:

\begin{Theorem}\label{thm:example}
  There exist a sequence of toric Mori fiber spaces with $m=n=2$ 
  such that $\mld(X)\to0$ and $\mld(Y)\approx C \cdot \mld(X)^4.$
\end{Theorem}

\begin{acknowledgments}
  This research originated in a discussion at the ACC workshop at the
  American Institute of Mathematics in May of 2012, as an attempt to
  answer a 
 question brought up by Yuri Prokhorov. The authors
  wish to thank Yu. Prokhorov, V.V. Shokurov
 and other participants of this workshop
  for many fruitful discussions.
\end{acknowledgments}

\section{Proofs of the main results}

We continue with the notation of the Introduction. We first need to
examine the combinatorics of the fans of $X$ and $Y$. We refer to
\cite{Fulton_ToricVarieties} or \cite{Oda_ConvexBodies} for the
general theory of toric varieties. We work over $\bC$ for simplicity,
although, as usual in toric geometry, the results remain true over a
field of positive characteristic as well. 

Recall that a toric variety $X$ is given by a pair $(N_X,\Sigma_X)$
where $N_X$ is a lattice (called the lattice of valuations) and
$\Sigma_X$ is a rational polyhedral fan in $N_X\otimes \bR$. A toric
map from a toric variety $X$ to a toric variety $Y$ is given by a
linear map $F\colon N_X\to N_Y$ such that its extension 
$\FR\colon N_X\otimes \bR\to N_Y\otimes \bR$ sends every cone in the
fan $\Sigma_X$ to
 a cone in the fan $\Sigma_Y$.

We denote by $N_Z$ the lattice $\Ker (F)$, and by $\Sigma_Z$ the
restriction of $\Sigma_X$ to $\Ker (F_{\bR})$. We recall the following
basic facts:

\begin{Fact}
  The morphism $f\colon X\to Y$ is proper iff $F_{\bR}\inv (\Supp\Sigma_Y)
  = \Supp\Sigma_X$.
\end{Fact}

\begin{Fact}
  A general fiber of $f\colon X\to Y$ is a product of a torus of
  dimension $\dim N_Z$ with the finite part, the product of finitely
  many copies of the group schemes $\mu_{r_i} = \Spec
  k[z]/(z^{r_i}-1)$. The character group of the finite part is the
  torsion subgroup of $\coker (F\colon N_X\to N_Y)$.
\end{Fact}

A Mori fiber space $f\colon X\to Y$ is a surjective proper morphism
with connected fibers, and a general fiber is connected and reduced.
Therefore, in our situation one has $\FR\inv( \Supp\Sigma_Y) =
\Supp\Sigma_X$, and the morphism of lattices $F\colon N_X\to N_Y$ is
surjective.

\begin{Fact}
  A toric variety $X$ is $\bQ$-Gorenstein, i.e. the canonical divisor
  $K_X$ is $\bQ$-Cartier iff there exists a function
  $\ell=\ell_{-K_X}\colon \Supp\Sigma_X\to \bR$ which is linear on each
  cone $\sigma\in\Sigma_X$ and such that
  $\ell(P_i)=1$ for each shortest integral generator $P_i$ of each ray
  $R_i$ of $\sigma$. 
\end{Fact}

\begin{Fact}
  A toric variety $X$ is $\bQ$-factorial, i.e. every Weil divisor is
  $\bQ$-Cartier iff the fan $\Sigma$ is simplicial, i.e. every cone is
  a simplex. 
\end{Fact}

\begin{Fact}\label{fact:mld}
  The $\mld(X)$ is computed as the minimum of the piece-wise linear function
  $\ell=\ell_{-K_X}$ on $\Supp\Sigma_X \cap N_X\setminus \{0\}$.
\end{Fact}

Obviously, our problem is local on $Y$, so we can assume that $Y$ is
affine. Since $X$ and $Y$ are $\bQ$-factorial, $Y$ is a quotient of
$\bA^n$ by a finite abelian group. Combinatorially, it is obtained from
the standard cone $C=\{(x_1,...,x_n | \forall x_i\geq 0\}$ in $\bR^n$,
with the lattice being a finite extension of the standard lattice
$\bZ^n$. Thus, the fan $\Sigma_Y$ consists of the cone $C$ and its faces.

If the shortest integral generators of the cone $C$ are the standard
basis vectors $e_i$ then the linear function computing $\mld(Y)$ is
simply $\ell_{-K_Y}=\sum x_i$.

\begin{Proposition} Suppose that $f:X \to Y$ is a toric Mori fiber space,
  with $\Q-$factorial $X$ and affine $Y$ as above. Denote by $F:N_X
  \to N_Y$ the map of the corresponding lattices, and extend it to the
  linear map $F_{\bR}$ from $N_X\otimes \bR = \bR ^{n+m}$ to
  $N_Y\otimes \bR = \bR^n.$ We choose the basis of $N_X\otimes \bR$ so
  that the map $F_{\bR}$ is the projection of $\bR^{n+m}$ to the last
  $n$ coordinates. Then the following is true about the fan of $X$ in
  $N_X\otimes \bR $.

  1) It has exactly $(n+m+1)$ one-dimensional cones (rays) $R_i$, $i=0,1,...,m+n$ of
  which $R_0, \dotsc, R_{m}$ are in $\Ker (F_{\bR})$ and for all
  $i=m+1,...,m+n$ the ray $R_{i}$ is mapped by $F_{\bR}$ to the ray
  of $C$, which is the span of the standard basis vector $e_{i-m}.$

  2) Any generators of $R_i$ for $i=0,...,m,$ form a simplex that
  contains $0.$

  3) The cones in the fan of $X$ are precisely the simplicial cones
  generated by $R_i$ for $i\in S\subset \{0,1,...,n+m\}$, where $S$ does
  not contain $\{0,1,...,m\}$.
\end{Proposition}

\begin{proof}

The condition $F_{\bR}\inv(C) = \Supp\Sigma_X$ implies that
\begin{enumerate}
\item The fan $\Sigma_Z$ is complete,
  i.e. $\Supp\Sigma_Z=N_Z\otimes\bR$. In particular, $\Sigma_X$ has
  at least $m+1$ rays in $\ker\FR$.
\item For each of the $n$ rays of $C$, there exists at least one ray
  of $\Sigma_X$ lying over it. 
\end{enumerate}

Recall that the  Picard group of $\bR-$Cartier divisors on  a toric variety $X$ is the quotient of the
space of piece-wise linear functions modulo the space of linear
functions on $\Sigma_X$.  Since both fans are simplicial (because
$X,Y$ are both $\bQ$-factorial) and full-dimensional, the relative
Picard number $\rho(X/Y)$ is the difference between the number of rays
of $X$ and $Y$ minus the relative dimension, $m$.

Therefore, $\Sigma_X$ has $m+n+1$ rays. Thus, $\Sigma_X$ has no other
rays other than the $(m+1)+n$ rays listed above, and over each ray of
$C$ there exists a \emph{unique} ray of $\Sigma_X$. This proves (1).

Finally, for this set of $n+m+1$ rays there is only
  one simplicial fan with support $F_{\bR}\inv(C)$: the one described
  in (3). This proves (2) and (3).
\end{proof}

We now choose a basis in $N_X\otimes \bR$ so that the last $n$
coordinate vectors are the primitive elements of $N_X$ on the rays
$R_i, i\geq m+1,$ denoted by $P_i.$  For
$i=0,...,m$ we also denote by $P_i$ the primitive elements of
$N_Z=\ker(N_X\to N_Y)$ on the rays~$R_i.$ 

By the above Proposition, the fan $\Sigma_X$ is isomorphic to the
Cartesian product 
$\Sigma_Y \times \Sigma_Z$. 
Since a general fiber of $X\to Y$ is connected, the map $F\colon
N_X\to N_Y$ is surjective. Therefore, one has $N_X\simeq N_Y\times
N_Z$. However, one need \emph{not} have $(N_X,\Sigma_X) \simeq (N_Y,\Sigma_Y)
\times (N_Z,\Sigma_Z)$. 
In particular, it is possible that $F(P_i)$ are not
primitive in the lattice $N_Y.$

Denote by $\Delta$ the simplex with vertices $P_i$
in $\Ker (F).$ This structure defines the toric Fano variety $Z$ of Picard
number one, which is the generic fiber of $f.$ We choose the
coordinates in $N_X\otimes \bR =\bR^{n+m}$ so that the lattice
generated by $P_i$ is the standard $\bZ^{m} \subset \bR^m = \Ker
(F_{\bR});$ the lattice $N_Z$ is a finite extension of it.

\medskip

We now describe our basic strategy for the proof of Theorem 1. Recall
that $\mld(Y)$ 
of a toric variety $Y$
is computed as the minimum of the linear function $\sum \limits_{i=1}^{n}
x_i$ over the non-zero points of $N_Y \cap C$. 

According to Fact~\ref{fact:mld}, the mld of
a toric singularity is the minimum of the log discrepancies of the
non-zero points of the corresponding cone, where the log discrepancy
of a point is the value on it of the linear function that equals $1$
on the rays of the cone. Suppose that for $Y$ this minimum is achieved
at some point $A$. We want to prove that if the log discrepancy of $A$
is very small, there must exist a point  in $N_X$, in one of the
cones
 of $\Sigma_X$, for which the log discrepancy is also small (less than the
given $\eps$). To look for this  point, we take a preimage $P$  of $A$ in
$N_X$ (it is possible, because $N_X\to N_Y$ is surjective),
and consider its multiples
$P, 2P, ... , tP$ modulo the lattice $\bZ ^{n+m}$, for some $t$ to be
specified later. If the log discrepancy of $A$ is really small, then
we can choose a fairly large $t$ such that for all these points the
sum of the last $n$ coordinates is still small. By Dirichlet Box
Principle, we can choose two of these points to be close to each
other, and take their difference. If we subtract in the correct order,
this produces a point $Q$ in $N_X$ with the last $n$ coordinates
nonnegative and with small sum; and the projection to the first $m$
coordinates being near the origin. Because the union of the cones for
$Z$ is the whole $\bR ^m,$ this projection must belong to some cone,
which implies that $Q$ lies in some cone for $X$ and has a small log
discrepancy there.

\medskip

To illustrate the method, we first establish Theorem 1 in the
particular case when $Z$ is an unramified in codimension one quotient of the usual
projective space $\bP^m$.
This means that the
barycentric coordinates of $0$ in the simplex $\Delta$ are
$(\frac{1}{m+1} , ... , \frac{1}{m+1}).$

\begin{Proposition} \label{prop:special-simplex} In the above
  notation, suppose additionally that the points $P_i$ for $1\leq i\leq m$
  are the standard $e_i \in \bR^{n+m}$, and $P_{0}=(-1,...-1;0,...,0)$
  (Here the semicolon separates the first $m$ coordinates from the
  last $n$). Then for any $\eps >0$, if $\mld(X)>\eps,$ then $\mld(Y) >
  \delta = (\frac{\eps}{2m})^{m+1}.$
\end{Proposition}
\begin{proof}

  Suppose that $\mld(Y)\leq \delta.$ Denote the point in $N_Y$ on which
  the mld is achieved, by $A.$ In other words, $A=(a_1,...a_n),$ where
  $a_i$ are nonnegative, not all zero, and $\sum a_i \leq \delta.$
  Because $F$ is surjective, $A=F(P)$ for some $P\in N_X.$ Suppose
  $P=(b_1,...,b_m; a_1,...a_n).$ We may additionally assume that all
  $b_i$ are in $[0,1),$ because $\bZ^m \subseteq N_Z.$

  Choose $t=\delta^{-\frac{m}{m+1}}$.  For all integers $k\in [0,t)$
  consider the points $P_k=kP \mod \bZ^{n+m} = (\lfloor kb_1 \rfloor
  ,..., \lfloor kb_m \rfloor; ka_1,...,ka_n)$ and their projections to
  $\bR^m$: $\bar{P_k}=(\lfloor kb_1 \rfloor ,..., \lfloor kb_m
  \rfloor)$.

  \begin{Lemma} Suppose for all integer $k\in [0,t)$,
    $Q_k=(b_{1,k},...,b_{m,k})$ are arbitrary points in
    $[0,1)^{m}$. Then there exist $i$ and $j$ so that for all
    $l=1,...,m$ we have $|b_{l,i}-b_{l,j}| \leq t^{-1/m}.$
  \end{Lemma}

\begin{proof}[Proof of Lemma]
  Identify $[0,1)^{m}$ with the quotient $\bR^m/\bZ^m,$ with the usual
  Haar probability measure. For each $Q_i$ consider a closed box
  neighborhood of it defined by the conditions $x_l\in
  [b_{l,i}-\frac{1}{2}t^{-1/m}, b_i+\frac{1}{2}t^{-1/m}] \mod Z.$ The
  volume (i.e. the Haar measure) of each such box is $t^{-1}.$ Note
  that the total number of points is $\lfloor t\rfloor > t,$ so the
  total sum of the volumes is greater than $1$. Thus there exist $i$
  and $j$ such that the corresponding boxes intersect. The triangle
  inequality in $\bR/\bZ$ implies the result.
\end{proof}

We apply the above Lemma to the points $Q_k=\bar{P_k}.$
Without loss of generality, we can assume that $i<j.$ Consider the
point $Q=P_j-P_i\in N_X.$ In coordinates, $Q=(\lfloor jb_1\rfloor-
\lfloor ib_1\rfloor , ..., \lfloor jb_m\rfloor -\lfloor ib_m\rfloor;
(j-i)a_1,..., (j-i)a_n)$. Note that $0< j-i \leq t,$ so the sum of the
last $n$ coordinates of $Q$ is at most $t\delta$. Suppose that for
$l\geq m+1$ we have $F(P_l) = c_l\cdot e_{l-m}$. Then the
contribution to the log discrepancy of $Q$ from the last $n$
coordinates is
$$k\sum_{l=m+1}^{n+m} \frac{a_{l-m}}{c_l} \leq k (a_l+...+a_n) \leq k\cdot \delta.$$

The first $m$ coordinates of $Q$ are less than $t^{-1/m}$ in absolute
value. Denote by $\bar{Q}$ the natural projection of $Q$ to $N_Z\otimes {\bR}$:
$$\bar{Q}=  (\lfloor
jb_1\rfloor- \lfloor ib_1\rfloor , ..., \lfloor jb_m\rfloor -\lfloor
ib_m\rfloor) = (q_1,...,q_m).$$
Then $\bar{Q}$ belongs to one of the cones of the fan for $Z$ as
follows. 

Case 1. All $q_l$ are nonnegative. Then $\bar{Q}$ belongs to the cone
$x_i\geq 0,$ which is the span of $P_i,$ for $1\leq i\leq m.$ The
contribution to the log discrepancy from the first $m$ coordinates is
at most $m\cdot t^{-1/m}$.

Case 2. At least one of the numbers $q_l$ is negative. Without loss of
generality, we can assume that $q_1$ is the smallest (i.e. the most
negative) of $q_l$. Then $\bar{Q}$ lies in the span of
$P_0;P_2,...,P_{m}.$ Its coordinates in that basis are
$(-q_1;q_2-q_1,...,q_m-q_1).$ The contribution to the log discrepancy
from the first $m$ coordinates is at most $(2m-1)\cdot t^{-1/m}$.

Putting it together, the log discrepancy of $Q$ is at most
$(2m-1)t^{-1/m} + t\delta$. Since we chose
$t=\delta^{-\frac{m}{m+1}}$, we get the log discrepancy of $Q$ to be
at most $2m\delta^{\frac{1}{m+1}}\leq \eps,$ which contradicts
$\mld(X)> \eps.$

This completes the proof of Proposition~\ref{prop:special-simplex}.
\end{proof}

\begin{Remark} One can improve the above estimate slightly by choosing
  $t$ to be a suitable constant times $\delta^{-\frac{1}{m+1}}$, and by
  a more ``projectively symmetric'' estimates for $\bar{Q}$. But it will still
  give the result of the form $\delta \ge \const(m) \cdot \eps ^{m+1}$,
  and would make the exposition considerably more muddled.
\end{Remark}

\begin{proof}[Proof of Theorems~\ref{thm:main} and \ref{thm:fiber-fixed}]
A slight generalization of the above argument yields
Theorem~\ref{thm:fiber-fixed}. Indeed, suppose $Z$ is an arbitrary
toric Fano variety of dimension 
$m$ with the corresponding simplex $\Delta,$ and suppose that the
barycentric coordinates of $0$ in $\Delta$ are $y_1,y_2,...,y_{m+1}.$
We can fix the vertices $P_i, i\leq m+1,$ in $\bZ^m.$ As before, we
can take $t$ to be $\delta ^{-\frac{1}{m+1}}$. We apply the same Lemma
(though one can get a somewhat better estimate by generalizing scaling
the boxes, keeping the same volume). As a result, the absolute values
of all coordinates of the point $\bar{Q}$ are again at most
$t^{-\frac{1}{m}}= \delta^{\frac{1}{m+1}}$. So for each of the $(m+1)$
linear functions corresponding to the $m-$dimensional cones of the fan
for $Z,$ the log discrepancy for $\bar{Q}$ will be bounded by constant
multiple of $\delta ^{\frac{1}{m+1}}$. The same estimate as above
proves that for the fixed $y_1,y_2,...,y_{m+1}$ one can choose
$\delta=\const\cdot \eps^{m+1},$ thus proving
Theorem~\ref{thm:fiber-fixed}. 

Finally, Theorem~\ref{thm:main} follows
from Theorem~\ref{thm:fiber-fixed} by a simple observation that if
$\mld(X) > \eps,$ then also $\mld (Z)>\eps$. 
By the main result of
\cite{BorisovBorisov} (BAB Conjecture for toric varieties) there are
only finitely many possible Fano varieties $Z$ with $\mld(Z)>~\eps$. 
\end{proof}

While the above argument may seem to imply the existence of a general
estimate for the $\mld(Y)$ in terms of $\mld(X) $ in the form $\const(m)
\cdot \eps^{m+1},$ the constant depends implicitly on $\eps$. In fact,
one simply cannot hope for the estimate above, in light of the
following example which proves Theorem 3.

\begin{Example} We fix $n=m=2$. Suppose $l$ is a natural
  number. Consider a triangle in $\bZ^2$ with vertices $(1,0),
  (-(l-1),1), (-(l-1), -1)$. This gives a weighted projective space;
  we multiply it by $\bA^2,$ and consider the quotient by the group
  $\mu_r,$ where $r=l^4+1,$ given by the weights
  $\frac{1}{r}(l,l^2;1,1).$ In other words, we take a lattice
  $\bZ^2\subset \bZ^4,$ and enlarge the latter by adjoining the point
  $\frac{1}{r}(l,l^2;1,1).$ The rays are $(1,0,0,0), (-(l-1),1,0,0),
  (-(l-1), -1, 0,0); (0,0,1,0), (0,0,0,1).$ The map $F$ is just the
  projection to the last two coordinates. The variety $Y$ is a cyclic
  quotient singularity of type $\frac{1}{r}(1,1).$
\end{Example}

We claim that for the above Example the $\mld(Y)$ is asymptotically
$\frac{1}{l^4},$ while $\mld(X)$ is asymptotically at least
$\frac{1}{2l}$. This would obviously imply Theorem 3.

The first part is easy: $\mld(Y)=2/r,$ which is asymptotically $2/l^4$.

For the estimate on $\mld(X)$, consider the point
$N=\frac{1}{r}(l,l^2;1,1)$ in $\bR^4$. We need to prove that no sums
$kN$ and points of $Z^4$ have small log discrepancy, in any of the
cones of $X$. Consider such point $Q=kP+B,$ where $k$ is an integer
from $1$ to $r-1$ and $B\in \bZ^4$. Clearly, we can assume that the
last two coordinates of $B$ are zero, thus $B\in \bZ^2\subset \bZ^4.$
Note the following.

1) If $k>(l^3)/2,$ then the contribution from the last two coordinates
is already too big. So we are only concerned with $k\leq (l^3)/2.$

2) Since $k\leq (l^3)/2,$ the first coordinate in $kN$ is between $0$
and $1/2.$ Therefore the points in the left cone are of no concern:
they would have log discrepancy contribution from the first two
coordinates at least $(1/ 2)/(l-1)$. For the points in the upper or
lower cone, if $k> l^2 /2,$ then the log discrepancy is at least $kl/r
> (l^3)/(2r)$, which is about $1/(2l)$. So we only need to consider $k
\leq l^2 / 2.$

3) Since $k \leq l^2 / 2,$ the second coordinate of $kN$ is between
$0$ and about $1/2.$ This rules out points in the lower cone. For the
upper cone, we clearly only need to be concerned with the points $kN+
(0,0)$.  And there the smallest $(x_1+l x_2)$ value is at least
$1\times l^2/r=l^3/r,$ which is about $1/l.$

\section{Miscellaneous remarks}

It may seem like one cannot avoid using the BAB conjecture to prove
Theorem 1. However, there is an explicit version of the toric BAB
theorem (proved by Lagarias and Ziegler \cite{LagariasZiegler}, and
originally by Hensley \cite{Hensley}, before \cite{BorisovBorisov})
which may probably be used to get an explicit bound of the form
$\delta=C(m)\eps^{d(m)}$. However this is by no means automatic, and
the correct power $d(m)$ is highly mysterious. Probably, for $m=2$ it
is $4$, but in higher dimensions the answer is not obvious.

By a more careful generalization of the argument for $\bP^m,$ one can
get an estimate for $\delta$ in terms of $m,$ $\eps$ and the Tian's alpha
invariant of $Z$ (that essentially measures ``asymmetry'' of the simplex
$\Delta$). Perhaps a generalization of this argument to non-toric case
will naturally use this invariant as well.

We also note a subsequent to our paper preprint \cite{Birkar}, whose
stated result is a proof of Shokurov's conjecture under the assumption
that the pair $(F,\Supp\Delta_F)$ belongs to a bounded family, where $F$
is a general fiber and $K_F+\Delta_F = (K_X+\Delta)|_F$, but without
the toric assumption.

\bibliographystyle{amsalpha} \renewcommand{\MR}[1]{}

\providecommand{\bysame}{\leavevmode\hbox to3em{\hrulefill}\thinspace}
\providecommand{\MR}{\relax\ifhmode\unskip\space\fi MR }
\providecommand{\MRhref}[2]{%
  \href{http://www.ams.org/mathscinet-getitem?mr=#1}{#2}
}
\providecommand{\href}[2]{#2}

\end{document}